
\def\quotation{%
  \par
  \begingroup
  \parindent=0pt
  \rightskip=2.3cm
  \leftskip=2.3cm
}
\def\endquotation{\par\endgroup}



\def\btc{%
  \leavevmode
  \vtop{\offinterlineskip 
    \setbox0=\hbox{B}%
    \setbox2=\hbox to\wd0{\hfil\hskip-.03em
    \vrule height .3ex width .15ex\hskip .08em
    \vrule height .3ex width .15ex\hfil}
    \vbox{\copy2\box0}\box2}}




\baselineskip=14pt
\parskip=10pt
\def\halmos{\hbox{\vrule height0.15cm width0.01cm\vbox{\hrule height
  0.01cm width0.2cm \vskip0.15cm \hrule height 0.01cm width0.2cm}\vrule
  height0.15cm width 0.01cm}}

\magnification=\magstephalf

\def\1{{\overline{1}}}
\def\2{{\overline{2}}}
\parindent=0pt
\overfullrule=0in

\def\frac#1#2{{#1 \over #2}}
\centerline
{\bf A Combinatorial-Probabilistic Analysis of Bitcoin Attacks}

\bigskip
\centerline
{\it Evangelos GEORGIADIS and Doron ZEILBERGER}
\bigskip

{\bf Abstract}: {\it In 2008, Satoshi Nakamoto}\footnote{*}{From here onwards, we refer to Satoshi Nakamoto as Dr.~Satoshi Nakamoto or Dr.~Nakamoto.} {\it famously invented bitcoin, and in his (or her, or their, or its) white paper
sketched an approximate formula  for the probability of a successful double spending attack by a dishonest party. 
This was corrected by Meni Rosenfeld, who, under more realistic assumptions, gave the exact probability (missing a foundational proof); and another formula (along with foundational proof), in terms of the Incomplete Beta function, was given later by Cyril Grunspan and Ricardo P\'erez-Marco, 
that enabled them to derive an asymptotic formula for that quantity.

Using Wilf-Zeilberger algorithmic proof theory, we continue in this vein and present a recurrence equation for the above-mentioned probability of 
success, that enables a very fast compilation of these probabilities.
We next use this recurrence to derive (in algorithmic fashion) higher-order asymptotic formulas, 
extending the formula of Grunspan and  P\'erez-Marco who only did the leading term. We 
then study the statistical properties (expectation, variance, etc.) of the duration of a successful attack.}

{\bf Important}: This article is accompanied by a Maple package, {\tt Bitcoin.txt}, available from the front of this article

{\tt http://sites.math.rutgers.edu/\~{}zeilberg/mamarim/mamarimhtml/bitcoin.html} \quad ,

where readers can also find  numerous input and output files.

In due course, the current package will also be available along with an implementation in MathCognify's own symbolic 
language\footnote{$^{\btc}$}{Donations in BTC are appreciated, if you'd like to support this type of research; 
some parts of the work will see the light of Free Software. The BTC address is: {\tt 3B69VRGSGt41K5GBrvCWaiUZ6uE2md9i9q}.} from

{\tt https://github.com/MathCognifyTechnologies/fr-crypto-bitcoin-1}. 

{\bf A Two-Phase Soccer Match}

In order to make this article self-contained and focus on the core combinatorial structure that underpins much of the process of the double spend attack, 
we will postpone the malevolent language of cyber-attacks until the penultimate section of this article, i.e., the notes and remarks section. 
There we provide context and references to the double-spending attack for the interested reader; additionally, 
we pinpoint weaknesses and inconsistencies in Dr. Satoshi Nakamoto's paper -- strengthening our belief that he is (or was) 
a competent scientist but not a combinatorialist. For now, we provide an equivalent model featuring a two-phase Soccer match, 
where one of the two teams is worse than the other.

There are two Soccer teams called the {\bf Good Androids}\footnote{$^\infty$}{Androids are humanoid robots. Our incentive for picking androids, 
is at least two-fold. For one, our creations are emotionless -- thereby strengthening some of our key assumptions below by taking the almost unpredictable 
nature or volatility of human emotions out of the parameter space. For another, networked androids enable easier generalizations.} {\bf Team} (henceforth {\bf G}) and the {\bf Bad Androids Team} (henceforth {\bf B}).
Luckily team G is better than team B. Whenever they play the probability of team B scoring the next goal is $q$,
with $q<\frac{1}{2}$, and hence the probability of team G scoring the next goal is $1-q$. The scoring of
any given goal is {\bf independent} of any past or future goals.

There is a positive integer $n$, decided beforehand. So the game has two parameters,  the continuous yet once chosen, fixed, $q$, with
$0<q<\frac{1}{2}$, and the discrete $n$. A forthcoming paper [E1], provides a more involved, protocol specific, analysis,
with a focus towards the stochastic version of $q$.

The match has two phases.

{\bf Phase I}: Play until team G scores $n$ goals. By then team B scored $m$ goals, say, and team G is $n-m$ goals ahead.

{\bf Phase II}: In the unlikely event that $m \geq n$, i.e., that team B scored at least as many goals as team G, team B is declared a winner
{\bf immediately}. Otherwise they keep playing until either team B caught up and tied the score, and is declared the winner,
or else team B is so much behind, say, by a zillion (not necessarily googol) goals, to make it hopeless for it to ever catch up, in which case it is
declared the loser.

{\bf Questions}: 

{\bf 1.} What is the probability, in terms of $q$ and $n$, of team B winning?

{\bf 2.} Assuming that team B won, what can you say about the random variable {\it duration of Phase II}?
In particular its {\bf expectation} and {\bf variance}?

In a beautiful, very lucid, article, Meni Rosenfeld ([R]), correcting an approximate formula in [N], 
answered Question 1 by stating and proving the following theorem.

{\bf Theorem 0} (Meni Rosenfeld, [R], p. 7, Eq. (1)): The probability of team B winning, henceforth
to be called the {\bf Rosenfeld Polynomial}\footnote{$^{**}$}{In honor of Meni Rosenfeld's pioneering spirit for correctly guessing that the negative binomial distribution provides a better approximation. Note, the duo, Cyril Grunspan and Ricardo P\'erez-Marco, are the first to cement this guess via proof in [GP,Proposition 5.3 on p.~10]; namely, that the negative binomial distribution is the exact distribution in the bitcoin model, under a given set of simplifying assumptions.}, $R_n(q)$, is given in terms of the following
expression, featuring {\bf binomial coefficients sums}:
$$
R_n(q) \, = \,
1 
\, - \, (1-q)^n \, \sum_{m=0}^{n-1} {{n+m-1} \choose {m}}\, q^m 
\, + \, q^n \, \sum_{m=0}^{n-1} {{n+m-1} \choose {m}}\, (1-q)^m  \quad .
$$

For the sake of completeness, let us engineer a proof for our model, phrased in terms of Soccer.

{\bf Proof of Theorem 0}: Let $m$ ($ m \geq 0$), be the number of goals scored by team  B at the end of Phase I.
The probability that it was indeed $m$ is given in terms of the {\bf negative binomial distribution}, and equals
$$
 {{n+m-1} \choose {m}}\,(1-q)^n \, q^m \quad .
$$
Indeed, the score right before (necessarily by team G) last goal was $n-1$ goals for team G and $m$ goals for team $B$.
There are  ${{n+m-1} \choose {m}}$ possible `histories' of getting there, since out of the total $n+m-1$ goals
scored, one has to {\bf choose} $m$ of them to go to Team B (and the remaining $n-1$ goals go to Team G).
By {\bf independence}, the probability of each such {\it history} is $(1-q)^n\,q^m$. The probability that
Team B won immediately at the end of Phase I is thus
$$
\sum_{m=n}^{\infty} {{n+m-1} \choose {m}}\,(1-q)^n \, q^m \, = \, 
1 \, - \, \sum_{m=0}^{n-1} {{n+m-1} \choose {m}}\,(1-q)^n \, q^m \, = \, 
1 \, - \, (1-q)^n \sum_{m=0}^{n-1} {{n+m-1} \choose {m}} \, q^m \quad .
$$

On the other hand if $m<n$, Team B still has a chance to catch up. Right now it is $n-m$ goals behind.
By the classical {\bf gambler's ruin} problem (e.g., [F] for a mathematically mature exposition or [p.~63, TB] 
as primer with worked out solution to problem 42, see also [Z1] for an algorithmic exegesis), 
the probability of catching up (i.e., tying the score) is $(\frac{q}{1-q})^{n-m}$. Hence if team B managed to score $m$ goals in Phase I, with $m<n$, the
probability of its succeeding to eventually catch up is:
$$
(\frac{q}{1-q})^{n-m} \, {{n+m-1} \choose {m}}\,(1-q)^n \, q^m  \, = \,
{{n+m-1} \choose {m}}\,(1-q)^{m}\, q^n \quad.
$$
Adding from $m=0$ to $m=n-1$ and combining with the probability of immediate success, done above, completes the
proof Theorem 0. $\halmos $

How would you compute the first, say, $10000$ Rosenfeld polynomials? An efficient way would be via
a {\it linear recurrence equation}, alias {\it difference equation}. Using the {\bf Zeilberger algorithm} [Z2], that
is part of {\bf Wilf-Zeilberger algorithmic proof theory} [PWZ], one  gets the following theorem.

{\bf Theorem 1}: The Rosenfeld polynomials $R_n(q)$ that compute the probability of team B winning if its
probability of scoring a single goal is $q$ and team G scored  $n$ goals in Phase I, satisfies
the following second-order linear recurrence with polynomial coefficients
$$
R_{{n}} \left( q \right) =-{\frac { \left( 4\,n{q}^{2}-4\,qn-6\,{q}^{2}
-n+6\,q+1 \right) R_{{n-1}} \left( q \right) }{n-1}}+2\,{\frac {q
 \left( q-1 \right)  \left( 2\,n-3 \right) R_{{n-2}} \left( q \right) }
{n-1}} \quad ,
$$
with initial conditions $R_0(q)=1$, $R_1(q)=2\,q$.

By finding an alternative expression for $R_n(q)$, in terms of the Incomplete Beta function\footnote{$^{@}$}{For completeness, we provide their formula in [GP, Theorem 6.1 p.~13], $P(z)=I_s(z,1/2).$},  
Cyril Grunspan and Ricardo P\'erez-Marco [GP] proved that $R_n(q)$ is asymptotically
$$
\frac{1}{\sqrt{\pi} \, (1-2q)} \frac{(4q(1-q))^n}{\sqrt{n}} \quad .
$$
Using the methods of  [Z3] (and the Maple package {\tt AsyRec.txt} available from there), we get a much more
precise asymptotics, with error term $O(\frac{1}{n^{\frac{11}{2}}})$ (and could easily go to much higher orders).
This is given by the following theorem.

{\bf Theorem 2}: $R_n(q)$ is asymptotically

$$
\frac {( 4\,q ( 1-q )  ) ^{n}}{\sqrt{\pi}\, (1-2q) } \, \cdot \, \frac{1}{\sqrt{n}} \, \cdot
$$
$$
( 1 \, + \, {\frac {12\,{q}^{2}-12\,q-1}{ 8 \, (  2\,q-1 ) ^{2} }} \, \cdot \, \frac{1}{n} \, + \,
{\frac {1}{128}}\,{\frac {400\,{q}^{4}-800\,{q}^{3}+120\,{q}^{2}+280\,q+1}{ ( 2\,q-1 ) ^{4} }} \, \cdot \, \frac{1}{n^2}
$$
$$
\, + \, {\frac {5}{1024}}\,{\frac {1344\,{q}^{6}-4032\,{q}^{5}+112\,
{q}^{4}+6496\,{q}^{3}-3444\,{q}^{2}-476\,q+1}{ ( 2\,q-1 ) ^{
6} }} \, \cdot \, \frac{1}{n^3} \, + \,
$$
$$
{\frac {21}{32768}}\,{\frac {20224\,{q}^{8}-80896\,{q}^{7
}-88832\,{q}^{6}+549632\,{q}^{5}-514400\,{q}^{4}+18368\,{q}^{3}+92368\,
{q}^{2}+3536\,q-1}{ ( 2\,q-1 ) ^{8} }} \, \cdot \frac{1}{n^4})  
$$
$$
+O(\frac{1}{n^{11/2}} ) \quad .
$$

If team B won, how long should it take? The following theorem gives an analog of Theorem $0$ for the expected duration of Phase II.

{\bf Theorem 3}: Let $A_n(q)$ be the following binomial coefficients sum
$$
A_n(q) \, = \, \, \frac{q^n}{1-2q} \, \sum_{m=0}^{n-1} {{n+m-1} \choose {m}}\, (1-q)^m  \, (n-m)\quad ,
$$
then the expected duration of Phase II, in case team B won,  let's call it $E_n(q)$, is given by
$$
E_n(q) \, = \, \frac{A_n(q)}{R_n(q)} \quad .
$$

{\bf Proof:} According to a classical result on the Gambler's Ruin ([F]), if a random walker starts out $L$ units to the right of $0$, 
and goes left one unit with probability $q$ and
right one unit with probability $1-q$, and $q<\frac{1}{2}$ then {\bf if} it makes it to $0$, then the expected time is $\frac{L}{1-2q}$.
(Recall that the probability that it happens is $(\frac{q}{1-q})^L$). Since if $m>n$  the duration of Phase II is $0$, the numerator of
the conditional expectation of the duration is $A_n(q)$, obtained by inserting $\frac{n-m}{1-2q}$ into the second sum in Theorem $0$.
Since everything is conditioned on team $B$ winning, we have to divide by $R_n(q)$.

Like for $R_n(q)$, the Zeilberger algorithm can be used to find a second-order recurrence satisfied by $A_n(q)$.
Since we already know how to compute $R_n(q)$ fast, this enables the fast compilation of a table for $E_n(q)$.
This is accomplished by the next theorem.

{\bf Theorem 4}: $A_n(q)$ satisfies the following second-order recurrence
$$
A_{{n}} \left( q \right) =-{\frac { \left( 4\,n{q}^{2}-4\,qn-6\,{q}^{2}
-n+6\,q \right) A_{{n-1}} \left( q \right) }{n-1}}+2\,{\frac {q \left( 
q-1 \right)  \left( 2\,n-3 \right) A_{{n-2}} \left( q \right) }{n-2}} \quad,
$$
subject to the initial conditions
$$
A_{{1}} \left( q \right) ={\frac {q}{1-2\,q}} \quad , \quad
A_{{2}} \left( q \right) =2\,{\frac {{q}^{2} \left( 2-q \right) }{1-2\,
q}} \quad .
$$

Using {\tt AsyRec.txt} ([Z3]) once again, our computer derived an asymptotic formula for $A_n(q)$ that we do not state. Combined with
the asymptotic formula for $R_n(q)$ given in Theorem 2, it yields the following theorem that gives an asymptotic formula for
the expected duration if team B wins (i.e., in case of a successful attack), let's call it $E_n(q)$.

{\bf Theorem 5}: The (conditional) expected duration of Phase II (in case team B won), where team G scores $n$ goals in Phase I, $E_n(q)$, 
tends to $\frac{1-q}{(1-2q)^2}$ as $n \rightarrow \infty$. More precisely,
the asymptotic expansion to order $4$ is:
$$
 \frac{( 1-q)}{( 1-2\,q )^2} \cdot
$$
$$     
(1 - {\frac {1}
{ ( 1- 2\,q )^2 }} \, \cdot \, \frac{1}{n} - {\frac {6\,{q}^{2}-6\,q-1}{ ( 1- 2\,q
 ) ^{4}}}\, \cdot \, \frac{1}{n^2} -{\frac {36\,{q}^{4}-72\,{q}^{3}+12\,{q}^{2}+24\,
q+1}{ (1-  2\,q ) ^{6}}} \, \cdot \, \frac{1}{n^3}
$$
$$
-{\frac {216\,{q}^{6}-648\,{q}
^{5}+276\,{q}^{4}+528\,{q}^{3}-306\,{q}^{2}-66\,q-1}{ ( 1- 2\,q
 ) ^{8} }} \, \cdot \, \frac{1}{n^4} )  
+ O(\frac{1}{n^5}) \quad .
$$

Using the same line of reasonings, one can derive a binomial coefficients sum for {\it any} given moment, and from it derive asymptotic formulas.
Doing this for the second moment, and combining with the asymptotic expression for the expectation (Theorem 5), produced the following theorem.

{\bf Theorem 6}: The variance of the (conditional) duration of Phase II (if team B won) tends to
$$
\frac{(1-q)(2-3q-4q^2)}{(1-2q)^4} \quad,
$$
as $n$ goes to infinity. More precisely, its asymptotic expansion to order $4$ is
$$
{\frac { ( q-1 )  ( 4
\,{q}^{2}-3\,q-2 ) }{ ( 2\,q-1 ) ^{4}}}-{\frac {
( q-1 )  ( 4\,{q}^{2}+2\,q-7 ) }{ ( 2\,q-1
 ) ^{6}}} \cdot \frac{1}{n} \,
- \,
{\frac { ( q-1 )  ( 24\,{q}^{4}-12\,{q
}^{3}-58\,{q}^{2}+29\,q+18 ) }{ ( 2\,q-1 ) ^{8}
}}\cdot \frac{1}{n^2} \,
$$
$$
-{\frac { ( q-1 )  ( 144\,{q}^{6}-216\,{q}^{5}-348\,{
q}^{4}+444\,{q}^{3}+328\,{q}^{2}-312\,q-41 ) }{ ( 2\,q-1
 ) ^{10}}} \cdot \frac{1}{n^3}
$$
$$
-{\frac { ( q-1 )  ( 864\,{q}^{8
}-2160\,{q}^{7}-1704\,{q}^{6}+4956\,{q}^{5}+4632\,{q}^{4}-9972\,{q}^{3}
+1586\,{q}^{2}+1711\,q+88 ) }{ ( 2\,q-1 ) ^{12}
}}\cdot \frac{1}{n^4}
 \, + \, O(\frac{1}{n^5}) \quad .
$$

Note that the standard deviation is {\bf larger}  than the expectation, hence, like the duration in the Gambler's Ruin problem, there is {\bf no } concentration about the mean.
\medskip

{\bf Notes and remarks}:

{\it Dr. Nakamoto's inconsistency}

For a survey of different types of attacks, see [CKLR]. We are interested in the double-spending attack as proposed by 
Dr.~Satoshi Nakamoto in [N]. Interestingly, or for that matter, amusingly, Dr.~Nakamoto's computational results are 
{\bf not consistent} with his explicit description. Dr.~Nakamoto states on page 7 in his white paper [N]:
\medskip
\quotation
The recipient waits until the transaction has been added to a block {\bf and $z$ blocks have been linked after it.}
\endquotation
\medskip

However, in his computation, if $z=0$, (i.e., the transaction in question is in a block, however, there are no subsequent blocks; also known as, 
1 confirmation) then, the probability of successfully double-spending is (miraculously !) 100\% 
(assuming that the attacker has a hashing power of a meager $q=0.1$ (fixed), that is, 10\%); refer to page 8 in [N].

The (pragmatic) validity of this statement seems to be a few standard deviations away from the mean, assuming a Gaussian distribution. 
In other words, this is far from the truth: when attempting to double spend a transaction with 1 confirmation, the attacker is surely not guaranteed success.

This inconsistency could arise from either our lack of truly understanding all implicit assumptions that Dr.~Nakamoto made, 
or it could be a result from an honest slip-up in formulating the explicit 
assumptions.\footnote{$^1$}{\it Other slip-ups of different nature, did occur in the paper; e.g., the omission of Feller's volume number 
in the reference section on page 9 in [N]. Any serious student of probability is aware that Feller published two volumes. 
Incidentally, the publication year that he provides, maybe indicative of his age or it's just noise. 
That said, in today's era, the printing year provided, for supposedly volume 1, is, harder to obtain.}

{\it Two (missing) implicit assumptions}

For what it's worth, our two satoshis, regarding implicit assumptions. There are at least two implicit assumptions that need to be made:
{\it
Not only does the attacker have to pre-mine a block (similar to the Finney attack as outlined by Hal Finney 
in the original post [HF]), but he is assumed to win propagation races.}

We note that this inconsistency isn't addressed by [GP] who exhibit mastery of conventional probabilistic machinery 
to attain their results.\footnote{$^3$}{Care needs to be taken when using technical jargon. In [GP] the term {\bf block validation} 
should be replaced by {\bf finding a solution to the proof-of-work} or an analogous phrase, since the 
process of {\bf validating blocks} is more or less instantaneous.}
Rosenfeld in [R], diplomatically circumvents this issue by explicitly outlining his assumptions.

At last, despite the inherent theoretical limitations of consensus in distributed computing as illustrated in [FLP], 
bitcoin's `consensus' computes at a market capitalization of $>$100 billion USD.\footnote{$^4$}{August 12,2018. {\tt www.coinmarketcap.com}}

{\bf Conclusion}: Using Wilf-Zeilberger algorithmic proof theory we investigated the probability, and duration, of successful 
bitcoin attacks. Readers can extend the above results
to higher moments using (for the time being) the Maple package {\tt Bitcoin.txt}, that also 
contains simulation programs that confirm all the theoretical formulas.
\bigskip

{\bf References}

[CKLR] M.~Conti, S.K.~E., L.~Lal, S.~Ruj,  {\it A survey on security and privacy issues of bitcoin}, in IEEE Communications Surveys \& Tutorials.
{\tt doi:10.1109/COMST.2018.2842460}, 2018.\hfill\break Available from:\hfill\break
{\tt http://ieeexplore.ieee.org/stamp/stamp.jsp?tp=\&arnumber=8369416\&isnumber=5451756} \quad . \hfill\break
[Accessed Aug. 15, 2018].

[E1] Evangelos Georgiadis, {\it Bitcoin: a probabilistic analysis I}. Forthcoming, 2018.

[F] W.~S.~Feller, {\it An introduction to probability theory and its applications: volume 1}, 3rd Edition (1968). New York: John Wiley \& Sons Inc.

[FLP] M.~J.~Fischer, N.~A.~Lynch, M.~S.~Paterson, {\it Impossibility of distributed consensus with one faulty process}, Journal of the ACM, 32(2):374--382, April 1985.

[HF] H.~T.~Finney II, {\it Re:Best practice for fast transaction acceptance - how high is the risk?}. {\tt https://bitcointalk.org/index.php?topic=3441.msg48384\#msg48384}, 11 Feb 2011.\hfill\quad\break [Accessed Aug. 15, 2018].

[GP] Cyril Grunspan and Ricardo P\'erez-Marco, {\it Double spend  races}, {\tt arXiv:1702.028672v2} [cs. CR], 17 Feb 2017. {\tt https://arxiv.org/abs/1702.02867}. \hfill\quad

[N] Satoshi Nakamoto, {\it Bitcoin: a peer-to-peer electronic system}, \hfill\break
{\tt https://bitcoin.org/bitcoin.pdf}, 2008. [Accessed Aug. 15, 2018].\quad

[PWZ] Marko Petkov{\v s}ek, Herbert S. Wilf, and Doron Zeilberger, {\it ``A=B"}, A.K. Peters, 1996. \hfill\break Available from:
{\tt https://www.math.upenn.edu/\~{}wilf/Downld.html} ~[Accessed Aug. 15, 2018].\hfill

[R] Meni Rosenfeld, {\it  Analysis of hashrate-based double spending}, {\tt arXiv:1402.2009v1} [cs.CR], 9 Feb 2014. {\tt https://arxiv.org/abs/1402.2009}. \hfill

[TB] D.~P.~Bertsekas and J.~N.~Tsitsiklis. {\it Introduction to Probability}, Athena Scientific, Nashua, NH, USA, second edition, 2008.\hfill

[Z1] Doron Zeilberger, {\it Symbol-crunching with the Gambler's Ruin problem},
in: ``Tapas in Experimental Mathematics", Contemporary Mathematics {\bf 457} (2008), 285-292, Tewodros Amdeberhan and Victor Moll, eds. 
\hfill\break 
Available from:
{\tt http://sites.math.rutgers.edu/\~{}zeilberg/mamarim/mamarimhtml/ruin.html}\hfill[Accessed Aug. 15, 2018].\quad

[Z2] Doron Zeilberger, {\it The method of creative telescoping}, J. Symbolic Computation {\bf 11}, 195-204 (1991). 
\hfill\break Available from:
{\tt http://sites.math.rutgers.edu/\~{}zeilberg/mamarimY/creative.pdf}\hfill\break[Accessed Aug. 15, 2018].\quad

[Z3] Doron Zeilberger, {\it AsyRec: A Maple package for Computing the Asymptotics of Solutions of Linear Recurrence Equations with Polynomial Coefficients},
The Personal Journal of Shalosh B. Ekhad and Doron Zeilberger, 2008. \hfill\break
{\tt http://sites.math.rutgers.edu/\~{}zeilberg/mamarim/mamarimhtml/asy.html}\hfill\break [Accessed Aug. 15, 2018].

\bigskip
\hrule
\bigskip
Evangelos Georgiadis, MathCognify Technologies Limited, Cheung Kong Center 19/F,
2 Queen's Road, Central, Hong Kong. \hfill\break
Email: {\tt egeorg at mathcognify dot com}   \quad 
\bigskip
Doron Zeilberger, Department of Mathematics, Rutgers University (New Brunswick), Hill Center-Busch Campus, 110 Frelinghuysen
Rd., Piscataway, NJ 08854-8019, USA. \hfill\break
Email: {\tt DoronZeil at gmail  dot com}   \quad 
\bigskip
\hrule
\bigskip
Aug. 15, 2018
\smallskip
(Revised Oct. 4, 2018)

\end